\documentclass[12pt]{amsart} 
\usepackage{latexsym}
\usepackage{amssymb} 
\usepackage{amsmath}
\usepackage{latexsym}
\usepackage{delarray}

 \newtheorem{thm}{Theorem}[section]

 \theoremstyle{definition}
 
 \theoremstyle{remark}

 \numberwithin{equation}{section}

  \newcommand{\Rn}{{{\mathbb R}^n}}

  \newcommand{\R}{{\mathbb R}}
  \newcommand{\Rone}{{\mathbb R}}

  \newcommand{\Cpx}{{\mathbb C}}
  \newcommand{\pa}{{\partial}}
  \newcommand{\rank}{{\; \rm rank \;}}

  \newcommand{\Imp}{{\;{\rm Im}\;}}
  \newcommand{\Rep}{{\;{\rm Re}\;}}

  \newcommand{\supp}{{\rm supp}}

\def\dslash{d\llap {\raisebox{.9ex}{$\scriptstyle-\!$}}}

\begin{document}
%
\title[On local and global regularity of Fourier integral operators]
 {On local and global regularity of Fourier integral operators}
\author{Michael Ruzhansky}

\address{%
Department of Mathematics\\
Imperial College London\\
United Kingdom}

\email{m.ruzhansky@imperial.ac.uk}

\thanks{The author was supported by a Royal Society grant
and by the EPSRC Grant EP/E062873/01}


\date{December 28, 2007}

\begin{abstract}
The aim of this paper is to give a review of local and global
properties of Fourier integral operators with real and complex
phases, in local $L^p$, global $L^2$, and in Colombeau's spaces.
\end{abstract}

\maketitle

\section{Introduction}

Fourier integral operators and their regularity properties
have been under study since their appearance in \cite{Ho71}.
One of the main motivations comes from the 
Cauchy problem for hyperbolic partial differential equations.
Consider a pseudo-differential operator
   \begin{equation*}
   P(t,x,\pa_t, \pa_x)=
   \pa_t^m+\sum_{j=1}^m P_j(t,x,\pa_x) \pa_t^{m-j}
   \label{eq:oper}
  \end{equation*}
 of order $m$. Here $t\in\Rone$ and $x\in X$, where 
 $X\subset\Rn$ is an open set.
 Operators $P_j(t,\cdot,\cdot)\in \Psi^j(X)$ are
 assumed to be 
 (classical) pseudo-differential operators of order $j$ on $X$.
 
 Let $\sigma_P(t,\tau,x,\xi)$ be the 
 principal symbol of $P$, i.e.
 the top order
 part of the full symbol of $P$, positively
 homogeneous of degree $m$ in $\xi$.
 As usual, $\tau$ and $\xi$ correspond to $D_t=-i\partial_t$ and
 $D_x=-i\partial_x$, respectively.

 Let $P$ be strictly hyperbolic, i.e. assume that 
 $\sigma_P$ can be factored as
 $$ \sigma_P(t,x,\tau,\xi)=\prod_{j=1}^m (\tau-\tau_j(t,x,\xi)),$$
 where $\tau_j(t,x,\xi)$ are real and distinct for $\xi\not=0$.
 Then they are positively homogeneous of degree one in $\xi$
 and smooth in $(t,x)$, for $\xi\not=0$.
 It is well-known that under 
 the strict hyperbolicity condition
 the Cauchy problem for $P$ is well-posed.
 The Cauchy problem for $P$ is the equation
  \begin{equation}\label{EQ:CP}
     \left\{ \begin{array}{ll}
         Pu(t,x)=0, & t\not=0, \\
         \pa_t^ju|_{t=0}=f_j(x), & 0\leq j\leq m-1.
      \end{array} \right.
  \end{equation}

 The loss of regularity for solutions $u(t,\cdot)$ compared
 to the Cauchy data depends on
 the operator $P$ and on the function spaces.
 For example, the usual energy conservation
 that holds in $L^2$ fails in $L^p$ for $p\not=2$.

 The general question that one is interested in is as
 follows. Let $W_j, W$ be some function spaces.
 The questions is whether for (small) $t$ we have the
 property that
 {$$ f_j\in W_j \textrm{ (for all } j \textrm{) }
 \Longrightarrow u(t,\cdot)\in W \;\;\; ?$$}
 First, recall the following examples of such mapping properties:
 \begin{equation}\label{EQ:mapping1}
  f_j\in L_{n-j}^{{ 2}} \Longrightarrow u(t,\cdot)\in 
    L^{{ 2}}_{n}
 \end{equation}
 and
 \begin{equation}\label{EQ:mapping2}
  f_j\in L^{{ p}}_{n-j+{ (n-1)|1/p-1/2|}} 
 \Longrightarrow u(t,\cdot)\in L^{{ p}}_n.
 \end{equation} 
 Here and in the sequel $L^p_s$ will denote the local Sobolev
 space over $L^p$ defined by the property that $f\in L^p_s$
 if $(1-\Delta)^{s/2}f\in L^p$.
 
 Property \eqref{EQ:mapping1} is the usual conservation of
 energy for hyperbolic equations in $L^2$. Property
 \eqref{EQ:mapping2} shows that the loss of 
 regularity occurs in $L^p$ for $p\not=2$.
 In the case of the wave equation
 this was established in \cite{St71}, \cite{Mi80}, 
 \cite{Pe80}. In general, this property was proved in
 \cite{SSS91} (see also \cite{So93} and \cite{St93}). 
 We note that this is quite different
 from $L^p$--$L^{p^\prime}$ properties,
 for some overview of which see e.g. \cite{RS05}.

 The main idea for deriving estimates \eqref{EQ:mapping2}
 is to use the fact that solution $u(t,x)$ to the Cauchy
 problem \eqref{EQ:CP} can be written
 as a sum of Fourier integral operators 
 (see e.g. \cite{La57}, \cite{Du96}, \cite{Ho85},
 \cite{Eg93}). Thus, the study of
 the regularity properties of solutions to the Cauchy problem
 \eqref{EQ:CP} are reduced to the analysis of the regularity
 properties of Fourier integral operators. A more complicated
 reduction is also possible for certain classes of hyperbolic
 systems, even with variable multiplicities, for example
 for systems with microlocally diagonalisable principal
 part (see \cite{Ro80}, and \cite{KR07} for generic 
 conditions on characteristics).
 
 \section{Fourier integral operators with real phases}
 
  Let $X, Y$ be open sets in $\Rn$. One defines the class of
  Fourier integral operators $T\in I^\mu(X\times Y,C^\prime)$ by
  the (microlocal) formula
   \begin{equation}\label{EQ:FIO}
     Tf(x)=\int_Y\int_\Rn
   e^{i\Phi(x,y,\theta)} a(x,y,\theta) f(y)d\theta\;dy,
   \end{equation} 
   where amplitude $a$ is a smooth function satisfying
   $$a\in S^\mu: |\partial_{x,y}^{\alpha}\partial_\theta^\beta 
   a(x,y,\theta)|
    \leq C_{\alpha\beta}(1+|\theta|)^{\mu-|\beta|},$$
   for all multi-indices $\alpha, \beta$, and 
   $\Phi$ is a phase function.
   First, we assume that the phase is  
   real-valued and satisfies the following properties:
   {\begin{itemize}
    \item[1.] $\Phi(x,y,\lambda\theta)=
      \lambda\Phi(x,y,\theta),$ for all $\lambda>0$;
    \item[2.] $d\Phi\not=0$;
    \item[3.] $d_\theta\Phi=0$ is smooth 
      (e.g. $d_\theta\Phi=0$ implies
      $d_{(x,y,\theta)}\frac{\partial\Phi}{\partial\theta_j}$ are
      independent).
   \end{itemize}}

  It turns out that Fourier integral operators have 
  a useful invariant quantity, the canonical relation,
  which is a conic Lagrangian manifold in 
  $(T^*X\backslash 0)\times (T^*Y\backslash 0)$, equipped
  with the symplectic form $\sigma_X\oplus -\sigma_Y$,
  where $\sigma_X$ and $\sigma_Y$ are canonical symplectic
  forms in $T^*X$ and $T^*Y$, respectively.
  The canonical relation of a Fourier integral operator 
  $T$ is defined by
  $C^\prime=\{(x,\xi,y,-\eta): (x,\xi,y,\eta)\in C\}$, where
   $$ C=C_\Phi=\{(x,d_x\Phi,y,d_y\Phi): d_\theta\Phi=0\}$$
  is
  the wave front set of the integral kernel of the operator $T$.
  We have the canonical projections:
    \begin{equation*}
    \begin{array}{ccccc}
       T^*X & \stackrel{\pi_X}{\longleftarrow} &
        C \subset T^*X\times T^*Y& 
        \stackrel{\pi_Y}{\longrightarrow} & T^*Y. \\
      & & \Big\downarrow\vcenter{%
         \rlap{$\pi_{X\times Y}$}}
      & & \\
      & & X\times Y & &
    \end{array}
   \label{eq:in-defpr}
  \end{equation*}

   Operators appearing in the analysis of hyperbolic equations
   satisfy the so-called local graph condition. It means
   that projections 
   $$\pi_X:C\to T^*X\backslash 0,\quad 
     \pi_Y:C\to T^*Y\backslash 0$$ 
   are local symplectic
   diffeomorphisms. In fact, because of the symplectic
   structure, one of them is diffeomorphic if and only if the
   other one is. In terms of the phase function $\Phi$, it
   means that we have
   $$
   \det\partial_y\partial_\theta\Phi\not=0,$$
   or similarly $\det\partial_x\partial_\theta\Phi\not=0.$
   Thus, from now on we will assume that Fourier integral 
   operators that we consider are non-degenerate, i.e. that
   the local graph condition is satisfied.
   Note that the degenerate case is quite more subtle, see
   e.g. \cite{Ph94} for some overview.
 
   The local $L^2$--boundedness of non-degenerate
   Fourier integral operators of order zero
   is well-known (\cite{Es70}, \cite{Ho71}).
   In $L^p$--spaces with $p\not=2$ there is a loss of
   regularity depending on the index $p$. 
   Indeed, it was shown in \cite{SSS91} that
   non-degenerate Fourier integral operators of
   order $\mu$ are locally
   bounded in $L^p$, $1<p<\infty$, provided that
   $\mu\leq -(n-1)|1/p-1/2|$. This follows by interpolation
   between the local $L^2$--boundedness of Fourier integral 
   operators of order zero and local boundedness of
   Fourier integral operators of order $-(n-1)/2$ from the
   Hardy space $H^1$ to $L^1$, established in \cite{SSS91}.
   Alternatively, it was also established in \cite{Ta04}
   that non-degenerate Fourier integral operators of
   order $-(n-1)/2$ are of weak (1,1)--type.
   In Section \ref{SEC:Lp} we will discuss the local
   $L^p$--properties further.
  
 \section{Fourier integral operators with complex phases}
  
  Let us now briefly describe Fourier integral operators with
  complex phase functions. Again, let $X, Y$ be open subsets 
  of $\Rn$. Then a Fourier integral operator 
  of order $\mu$ with complex phase
   is defined by the same formulae as in
  \eqref{EQ:FIO}, namely
  \begin{equation}\label{EQ:FIOc}
   Tf(x)=\int_\Rn\int_\Rn
   e^{{\rm i}\Phi(x,y,\theta)} a(x,y,\theta) f(y)d\theta\;dy,
  \end{equation} 
 where amplitude $a\in S^\mu$ is of order $\mu$ and the 
 complex-valued phase function 
 $\Phi\in\Cpx$ satisfies conditions:
    { \begin{itemize}
    \item[1.] $\Phi(x,y,\lambda\theta)=\lambda
    \Phi(x,y,\theta),$ for all $\lambda>0$;
    \item[2.] $d\Phi\not=0$;
    \item[3.] $d_\theta\Phi=0$ is smooth
      (e.g. $d_\theta\Phi=0$ implies
      $d\frac{\partial\Phi}{\partial\theta_j}$ are
      independent over $\Cpx$);
    \item[4.] $\Imp\Phi\geq 0.$
    \end{itemize}}
 The last property is clearly necessary in general, for the
 integral \eqref{EQ:FIOc} to be well-defined.
 The theory of Fourier integral operators with complex phases
 was developed by Melin, Sj\"{o}strand in \cite{MS75} and
 \cite{MS76}, see also \cite{Ho85} and
 \cite{Tr82} for slightly alternative descriptions.
 These operators are also related to the Maslov canonical
 operator.

 In \cite{LSV94}, Laptev, Safarov and Vassiliev
 established global parameterisations of Fourier integral 
 operators with real phases by globally defined complex
 phases (see also \cite{SV96}).
  
  Again, the canonical relation of $T$ 
  with complex phase $\Phi$ is defined by
   $$ C=C_\Phi=\{(x,d_x\Phi,y,d_y\Phi): d_\theta\Phi=0\}.$$
  Now we can identify $C$ with a subset of
  $\widetilde{T^*(X\times Y)}$, where
  $\widetilde{T^*(X\times Y)}$ is
  an almost analytic extension of $T^*(X\times Y)$.
  Here one needs to work with almost analytic extensions of
  real sets, for example in order to interpret
  the eikonal equation $\partial_t\Phi=\tau(t,x,\nabla_x\Phi)$
  when $\Phi$ is complex-valued.

  Local $L^2$ properties for operators with complex phases 
  are more subtle. For example, in \cite{MS76},
  Melin and Sj\"{o}strand showed that operators of order zero
  are locally bounded in $L^2$ under a complex-valued version
  of a local graph condition. This condition and its extension
  will be discusses in the next section in more detail.
  A more general result was established by H\"{o}rmander 
  in \cite{Ho83}. Assume that  
  for every $\gamma\in C_\Rone$, the maps
   $(T_\gamma C)_\Rone$ to $T_{\pi_X(\gamma)}T^*X$ and
   to $T_{\pi_Y(\gamma)}T^*Y$ are injective. Then it is
   shown in \cite{Ho83} that Fourier integral operators
   of order zero with complex-valued phases are locally bounded in
   $L^2$. In fact, one can shows that if the mappings
   $$C_\Rone\to T^*X\backslash 0, 
    \quad C_\Rone\to T^*Y\backslash 0$$ are injective,
   then $T\circ T^*\in \Psi^0_{1/2,1/2}$ is a
   pseudo-differential operator of order zero and type
   $(1/2,1/2)$. This is 
   one of the differences with operators with real-valued phases,
   where $T\circ T^*\in \Psi^0_{1,0}$ 
   is a pseudo-differential operator
   of order zero and type $(1,0)$.
 
  \section{Local $L^p$ regularity properties of 
     Fourier integral operators}
  \label{SEC:Lp}

   The local $L^p$ properties of Fourier integral operators
   can be summarised in the following table:
       \begin{center}
          \begin{tabular}{|c||c|c|} \hline
             &  {$\Phi\in\mathbb R$} & {$\Phi\in\mathbb C$} \\ \hline
          $L^2$ & \cite{Es70,Ho71} & 
          \cite{MS75,MS76,Ho83} \\ \hline
          $L^p$ & \cite{SSS91} & \cite{Ru01} \\ \hline
           \end{tabular}
       \end{center}

   To explain this is more detail, we have the following 
   results for the local $L^p$--continuity:
    \begin{itemize}
      \item  Real $L^2$ theory: 
        Eskin \cite{Es70} and H\"ormander \cite{Ho71}:
        operators of order zero with real phases
        are locally continuous in $L^2$.
      \item  Real $L^p$ theory:
        Stein \cite{St71}, Beals \cite{Be82}, Miyachi \cite{Mi80},
        Peral \cite{Pe80}, Sugimoto \cite{Su92}, and
        finally Seeger, Sogge and Stein \cite{SSS91}:
        operators of order $\mu$ with real phases are 
        locally continuous in $L^p$, provided that 
        $\mu\leq -(n-1)|1/p-1/2|$.
      \item Complex $L^2$ theory:
        Melin and Sj\"ostrand \cite{MS75, MS76} and
        H\"ormander \cite{Ho83}:
        operators of order zero with complex phases 
        are locally continuous in $L^2$.
      \item Complex $L^p$ theory: Ruzhansky \cite{Ru01}:
        extension of the properties above to derive the 
        unified local $L^p$ properties for 
        operators  with complex phases.
    \end{itemize}
   We will now review local   
   $L^p$ results for Fourier integral operators with complex 
   phases in more detail.
  
  Let us introduce the following local graph condition for 
  complex-valued phases. We will assume that there
  exists some $\tau\in\Cpx$ such that we have
 $$ {\rm (LG)}\;\;\quad \quad 
   \det\partial_\theta\partial_y(\Rep\Phi+
   \tau\Imp\Phi)\not=0.$$
 Since the equation of the determinant to zero is polynomial
 in $\tau$, we can conclude that this condition 
 is equivalent to the existence of $\lambda\in\Rone$ such that
 $$\Psi=\Rep\Phi+\lambda\Imp\Phi$$ 
 defines a (real) local
  graph, i.e. that
  { $\det\partial_\theta\partial_y\Psi(x,y,\theta)\not=0.$} 
 Condition (LG) weakens Melin--Sj\"ostrand's local
 graph condition, which corresponds to the case $\tau=i$.
 Following \cite{Ru01}, under this condition we have the 
 following results. For the sake of simplicity in notation
 we will suppress primes in the notation for canonical relations.

\begin{thm}[\cite{Ru01}]
    Let $C\subset (T^*X\backslash 0)\times 
    (T^*Y\backslash\widetilde{0)}$ be
  a smooth complex positive homogeneous canonical relation,
  closed in $\widetilde{T^*(X\times Y)\backslash 0}$. 
  Assume that (LG) holds.
  Let $\mu=-(n-1)|1/p-1/2|$,
  $1<p<\infty$. Then
  $T\in I^\mu(X\times Y,C^\prime)$ is continuous
  from $L^p_{comp}(Y)$
  to $L^p_{loc}(X)$.
    \label{th:thfio1}
\end{thm}
 This theorem 
 extends both the result of
 Melin and Sj\"ostrand in \cite{MS75, MS76} ($L^2$, $\Phi\in\Cpx$) and
 the result of Seeger, Sogge and Stein in \cite{SSS91}
 ($L^p$, $\Phi\in\Rone$). 
 We can also note that Theorem \ref{th:thfio1} can not be 
 obtained from its real-valued version since
 we only have that
 there exists a real conic Lagrangian manifold
 $C_0\subset T^*(X\times Y)\backslash 0$ 
 such that
 $C_\Rone\subset C_0$ and such that the class of operators with
 complex phase of order $\mu$ with canonical relation $C$
 is contained in the class of
 operators of the same order $\mu$ with real phases
 with real canonical relation $C_0$ of type
 $(1/2,1/2)$, i.e. we have $I^\mu_\rho(X,Y;C)\subset
 I^\mu_{1/2}(X,Y;C_0)$. Thus, we may only conclude that 
 operators with complex phases of order $\mu$ are locally
 bounded in $L^p$ provided that $\mu\leq -(n-1/2)|1/p-1/2|$.
 This order is clearly worse 
 than the order $\mu$ given in
 in Theorem \ref{th:thfio1}.

 \section{Cauchy problem for operators with complex 
    characteristics}
  The results can be applied to establish the local regularity
  properties to solutions to the Cauchy problem for operators
  with complex characteristics.
  Consider a (classical) pseudo-differential operator
   \begin{equation*}
   P(t,x,D_t, D_x)=D_t^m+\sum_{j=1}^m P_j(t,x,D_x) D_t^{m-j}
   \label{eq:oper-c}
  \end{equation*}
 of order $m$, where $t\in [0,T]$ and $x\in X\subset\Rn$.
 Here $P_j\in C^\infty([0,T],\Psi_{cl}^j(X))$ are classical
 pseudo-differential operators of order $j$.

 The principal symbol of operator $P$ is given by
 $$ \sigma_P(t,x,\tau,\xi)=
   \tau^m+\sum_{j=1}^m P_j(t,x,\xi) \tau^{m-j}=
   \prod_{j=1}^m (\tau-\tau_j(t,x,\xi)).$$
 The Cauchy problem for operator $P$ is the equation
  \begin{equation*}
     \left\{ \begin{array}{ll}
         Pu(t,x)=0, & t\not=0, \\
         \pa_t^j u|_{t=0}=f_j(x), & 0\leq j\leq m-1.
      \end{array} \right.
    \label{eq:Cauchy}
  \end{equation*}
  It was shown in \cite{Tr82} that the propagator for
  this Cauchy problem can be expressed in terms of
  Fourier integral operators with complex phases, under
  the assumption that characteristic roots $\tau_j\in\Cpx$
  are complex-valued and satisfy the following assumptions
  (A1) and (A2):
 $$(A1)\quad\textrm{Simple characteristics:}\; \tau_i\not=\tau_j,
  i\not=j, \xi\not=0.$$
 Therefore, we have that
 $\tau_j\in C^\infty([0,T]\times(T^*X\backslash 0), \Cpx)$
 is smooth and positively homogeneous of order one in $\xi$.
 We also assume that
 $$(A2):\quad \Imp\tau_j(t,x,\xi)\geq 0\; \textrm{ in }\;
    [0,T]\times(T^*X\backslash 0).$$

  \begin{thm}[\cite{Ru01}] Let $P$ satisfy {\rm (A1)} and {\rm (A2)}.
    Let $1<p<\infty$. Let $f_j$ be compactly supported. Then
    we have locally:
    $$f_j\in L^p_{\alpha-j+{ (n-1)|1/p-1/2|}}\; (j=0,\ldots,m-1)\;
     \Longrightarrow
      u(t,\cdot)\in L^p_\alpha.$$
   Also,  there exists a constant 
   $C_T$ such that for all $t\in[0,T]$ we have the local estimate
   $$||u(t,\cdot)||_{L^p_\alpha}\leq C_T\sum_{j=0}^{m-1}
    ||f_j||_{L^p_{\alpha-j+{ (n-1)|1/p-1/2|}}}.$$
    \label{th:thCauchy}
  \end{thm}
  {The result is sharp in general, since the loss of 
  { $(n-1)|1/p-1/2|$} derivatives
  is sharp for wave type equations (where $\tau_j\in\Rone$).}
 
  We note that the regularity results for operators with
  complex-valued functions have other applications, see e.g.
  \cite{Ru01} for an application to the oblique derivative
  problem (see also \cite{Ru03} and \cite{Ru01}).
 
 \section{Sharpness of local $L^p$ properties}

    Let $T\in I^\mu(X,Y;C)$ be a Fourier integral operator with
    a real-valed phase function $\Phi(x,y,\theta)$. 
    The question for the sharpness of local $L^p$
    estimates is to find the largest $\mu$ for which
    operators of order $\mu$ are continuous from
    $L^p_{comp}$ to $L^p_{loc}$. Alternatively, by using the
    calculus, one can look for the
    smallest $\alpha$ such that operators of order zero
    are continuous from the Sobolev space $L^p_\alpha$ to $L^p$. 
    In this case we have $\mu=-\alpha$.

    Let
    $$k:=\max_{x,y,\theta} \rank \frac{\partial^2}
    {\partial\theta^2}\Phi(x,y,\theta).$$
    In particular, the meaning of $k$ is that 
    the dimension of the singular support of the
    integral kernel of $T$ is less or equal to $n+k$
    (see \cite{Ru99a}). For example, we have
    $k=n-1$ for the
    solutions of the wave equation and $k=0$ for 
    pseudo-differential operators. The following result
    was proved in \cite{Ru99}: 
  \begin{thm}[\cite{Ru99}]\label{Th:sharpness}
    If a Fourier integral operator $T$ of order zero
    with real phase is locally 
    bounded from $L^p_{\alpha}$ to $L^p$, then
    $\alpha\geq\alpha_p=k|1/p-1/2|$.
  \end{thm}
   There are the following important cases of this result: 
   
   \begin{itemize}
   \item $p=2$; here $\mu=0$.
   \item $1<p<\infty$ and $k=n-1$: 
   this is the case of solutions to the wave equations
   (\cite{Pe80}) and more general hyperbolic equations
   (\cite{SSS91});  here $\mu=-(n-1)|1/p-1/2|$. 
   \item $k=0$: this is the case of pseudo-differential 
   operators; here $\mu=0$.
   \end{itemize}
  
  The orders in Theorem \ref{Th:sharpness} indicate that positive
  results on the $L^p$ boundedness and the order $\mu$
  in Theorem \ref{th:thfio1} can be improved if we impose
  further conditions on the canonical relation $C$.
  
  Let us first review such positive results for operators
  with real-valued phases.   
 Define $$\Sigma=\{(x,y): \rank \frac{\partial^2}{\partial\theta^2}
 \Phi(x,y,\theta)=k \}.$$ If $(x,y)\in\Sigma$, the level set
 of $\nabla_\theta\Phi$ is a linear space of dimension $n-k$
 in
   $C_\Phi\subset T^*(X\times Y),$ corresponding to a 
   linear subspace 
   in the conormal bundle
   $$ \begin{array}{l} N^*\Sigma=\{(x,y,\xi,\eta)\in T^*(X\times Y):
   (x,y)\in\Sigma,  \\
   \qquad\qquad \xi(\delta x)+\eta(\delta y)=0, \forall
   (\delta x,\delta y)\in T_{(x,y)}\Sigma\}.\end{array}$$
   The inclusion $N^*\Sigma\subset C_\Phi$ is dense, and
   the conormal bundle $N^*\Sigma$ consists of affine
   fibers.

   In \cite{SSS91}, the following so-called smooth factorization 
   condition (SF) was introduced. Assume that the mapping
   $$(x,y)\mapsto (\theta-\textrm{ level set of }
     \nabla_\theta\Phi(x,y,\theta))
     :\Sigma\to \mathbb{G}_{n-k}(\Rn)$$
   is smoothly extendible from $\Sigma$ to $\pi_{X\times Y}(C_\Phi)$,
   the singular support of the integral kernel of
   operator $T$ with canonical relation $C_\Phi$. 
   Here $\mathbb{G}_{n-k}(\Rn)$ denoted the Grassmanian which
   is the collection of all ($n-k$)--dimensional linear
   subspaces of $\Rn$.
   
   We can note the the smooth factorization condition
   (SF) automatically holds in the case
   $k=0$ (pseudo-differential operators) and in the case
   $k=n-1$ (hyperbolic equations).
  
  \begin{thm}[\cite{SSS91}]
    Under the smooth factorization condition {\rm (SF)} 
    operators Fourier integral operators of order zero with
    real phases are locally continuous from
    $L^p_{\alpha_p}$ to $L^p$, provided that 
    $\alpha_p= k|1/p-1/2|$ and $1<p<\infty$.
  \end{thm}
  
  Thus, a question arises when and whether the smooth factorization 
  condition (SF) is satisfied. One answer is given by the
  following theorem:

  \begin{thm}[\cite{Ru99a}]\label{Th:translationinv}
   For translation invariant 
   Fourier integral operators with analytic phases, 
   the smooth factorization condition {\rm (SF)} is satisfied,
   provided that $n\leq 4$, or $k\leq 2$.
  \end{thm}
  As a consequence, one obtains a sharp list of 
  $L^p$--properties for translation invariant
  Fourier integral operators in $\Rn$, with $n\leq 4$, or $k\leq 2$.
  The analysis is based on the theory of affine fibrations
  developed in \cite{Ru00}.

  Theorem \ref{Th:translationinv} shows that in order for the
  smooth factorization condition to break, one needs to have
  the dimension at least $n\geq 5$. 
  In such case examples
  for the failure of (SF) have been constructed in
  \cite{Ru99} and in \cite{Ru00}. A typical example is the
  phase function defined by
    $$\Phi(x,y,\theta)=\langle x-y,\theta\rangle+ \theta_1\theta_2^2\theta_5^{-2}+
     (\theta_3\theta_5-\theta_2\theta_4)^2\theta_5^{-3}.$$
  If the operator is not translation invariant, one can construct
  simpler examples already in lower dimensions. For example,
  if $n=3$, the phase function   
    $$\Phi(x,y,\theta)=\langle x-y,\theta\rangle+ 
    \theta_3^{-1}(y_1\theta_1+y_2\theta_2)^2$$
  provides an example of the failure of (SF). However, in
  \cite{Ru02} it was shown that the loss of regularity for
  the corresponding Fourier integral operators
  is the same as for operators satisfying the smooth factorization
  condition (SF). 
 
  An analogue of the smooth factorization condition
  for complex phases was discussed in \cite{Ru01}. The
  complex-valued phase function $\Psi$ satisfies the complex
  smooth factorization condition ($\Cpx$SF) if there
  exists a real $|\tau|<1/\sqrt{3}$ such that: 
  \begin{itemize}
    \item $\Psi=\Rep\Phi+\tau\Imp\Phi$ satisfies (LG); 
    \item $\Psi$ satisfies real (SF) for some $k$.
  \end{itemize}  
    
  \begin{thm}[\cite{Ru01}]
    Let $C\subset (T^*X\backslash 0)\times 
    (T^*Y\backslash\widetilde{0)}$ be
  a smooth complex positive homogeneous canonical relation,
  closed in $\widetilde{T^*(X\times Y)\backslash 0}$. 
  Assume that the complex smooth factorization condition
  {\rm ($\Cpx$SF)} is satisfied. 
  Let $\mu=-k|1/p-1/2|$ and
  $1<p<\infty$. 
  Then Fourier integral operators of order $\mu$ with
  complex phases and canonical relation $C$ 
  are continuous
  from $L^p_{comp}(Y)$
  to $L^p_{loc}(X)$.  
    \label{th:thfio2}
\end{thm}
Theorem \ref{th:thfio2} deals with operators with symbols of
type $(1,0)$. If the amplitude of a Fourier integral operator
with the canonical relation as in Theorem \ref{th:thfio2}
is in the class $S_{\rho,1-\rho}^\mu$, $\rho\geq 1/2$,
then it is continuous from $L^p_{comp}(Y)$
to $L^p_{loc}(X)$ provided that $\mu=-(k-(n-k)(1-\rho))|1/p-1/2|$.
  
  There are close relations between the regularity properties
  of Fourier integral operators with complex phases
  to the symplectic and analytic
  geometries, in particular to singularities of
  the so-called affine fibrations. We refer to \cite{Ru00}
  for details.
  
 \section{Global $L^2$ estimates for Fourier integral operators}

 We now consider
 globally defined Fourier integral operators on $\Rn$ defined by
  \begin{equation}\label{EQ:globalT}
   Tu(x)=\int_\Rn\int_\Rn
   e^{i(x\cdot\xi+\phi(y,\xi))} a(x,y,\xi) u(y)d\xi\;dy,
  \end{equation}  
 where real-valued phase function 
 $\phi\in C^\infty(\Rn\times\Rn)$ will satisfy conditions 
 specified below. We note that in the case of 
 pseudo-differential operators we have 
 $$\phi(y,\xi)=-y\cdot\xi.$$

 Similar results will hold for the adjoint operator.
 In particular, this includes operators of the form
 $$Su(x)=\int_\Rn e^{i\phi(x,\xi)} a(x,\xi) 
 \widehat{u}(\xi) d\xi,$$
 which appear as propagators for some classes
 of hyperbolic equations.

 Global properties on $L^2(\Rn)$ of pseudo-differential operators
 are well-known (see e.g.
 Calderon--Vaillancourt \cite{CV71}, Cordes \cite{Co75},
 Coifman--Meyer \cite{CM78}, Childs \cite{Ch76},
 Sugimoto \cite{Su88}, Boulkhemair \cite{Bo95},
 and many other contributions).
    
 The case of pseudo-differential operators
 was also studied but much less is known
 (see Asada--Fujiwara \cite{AF78},
 Asada \cite{As81, As84}, Fujiwara \cite{Fu75},
 Kumano-go \cite{Ku76}, Boulkhemair \cite{Bo97}).
 A certain disadvantage of these results is that
 in all these papers the authors made an assumption
 that $\partial_\xi\partial_\xi\phi$ is globally
 bounded on $\Rn\times\Rn$, which fails in many important situations. 
    
    For example, in a typical application to smoothing problems
    (as in \cite{RS04}--\cite{RS06c}),
    the canonical transforms used for the changes of
    variables on the Fourier transform side have the phase 
    of the form
    $x\cdot\xi-y\cdot \psi(\xi)$, in which case we have 
    $\phi(y,\xi)=y\cdot \psi(\xi)$,
    where { $\psi$ is positively homogeneous of order one}. 
    But then
    $$\partial_\xi\partial_\xi\phi(y,\xi)=
     y\cdot\partial_\xi\partial_\xi\psi(\xi)$$ is 
    unbounded on $\Rn\times\Rn$.
    
    There are several questions that arise.
    For example, what are the minimal growth conditions
    on the phase and amplitude for operators 
    \eqref{EQ:globalT} to be globally bounded on $L^2(\Rn)$.
    Moreover, if we want to have weighted estimates in 
    $L^2(\Rn)$, or weighted estimates in Sobolev spaces
    over $L^2(\Rn)$, the question arises of what are 
    the minimal requirements for the global calculus
    of such operators.

 Let us now assume that the phase satisfies 
 the following conditions
 on $\supp a$: 
 
 \bigskip
  (C1) $\qquad|\det\partial_y\partial_\xi\phi(y,\xi)|
   \geq C>0,$ 
        $\quad\forall (y,\xi)\in\Rn\times\Rn;$ 
        
        \bigskip
  (C2) $\qquad|\partial_y^\alpha\partial_\xi\phi(y,\xi)|\leq C_\alpha$,
        $|\partial_y\partial_\xi^\beta\phi(y,\xi)|\leq C_\beta$ \\ 
   ${ }$ $\qquad\qquad\qquad\forall (y,\xi)\in\Rn\times\Rn$, 
   $\quad 1\leq |\alpha|,|\beta|\leq 2n+2.$

  Note that condition (C1) is a 
  global version of the local graph condition, which is 
  necessary even
  for local $L^2$--bounds for operators with
  amplitudes in $S_{1,0}^0$. 
  
  The importance of condition (C2) is that now we must take
  only mixed derivatives with respect to $y$ and $\xi$,
  so the phase functions for the canonical transforms
  in smoothing problems satisfy this condition.
  Indeed, if
  $\phi(y,\xi)=y\cdot \psi(\xi)$, where $\psi$ is homogeneous of 
  order one for large $\xi$
  and $|\det D\psi(\xi)|\geq C>0$, then condition 
  (C2) is satisfied for large frequencies.
  An additional argument is required for small frequencies
  and it can be found in \cite{RS06a}.

  It can be noted that conditions (C1) and (C2) are 
  considerably weaker than those appearing
  in the analysis of SG--pseudo-differential operators
  (see Cordes \cite{Co95}),
  SG--Fourier integral operators, 
  (see Coriasco \cite{Co99}), or for operator
  of Shubin type (see Boggiatto, Buzano and Rodino \cite{BBR96}).

  We consider first operators of the form
  $$ 
   Tu(x)=\int_\Rn\int_\Rn
   e^{i(x\cdot\xi+\phi(y,\xi))} a(x,\xi) u(y)d\xi\;dy.
   $$
  For such operators we have the following theorem

 \begin{thm}[\cite{RS06a}]
  Let $\phi(y,\xi)$ satisfy conditions {\rm (C1), (C2)}. 
  Let $a(x,\xi)$
  satisfy one of the following conditions:

  (1)[Calder\'on-Vaillancourt]
    $\partial_x^\alpha\partial_\xi^\beta a(x,\xi)\in
      L^\infty(\R^n_x\times\R^n_\xi)$, 
      $\alpha,\beta\in\{0,1\}^n.$

  (2)[Cordes] $\partial_x^\alpha\partial_\xi^\beta a(x,\xi)\in
      L^\infty(\R^n_x\times\R^n_\xi)$, $|\alpha|,|\beta|\leq [n/2]+1.$

  (3)[Cordes] $\exists\lambda,\lambda^\prime>n/2:
    (1-\Delta_x)^{\lambda/2}(1-\Delta_\xi)^{\lambda^\prime/2}a(x,\xi)
     \in L^\infty(\R^n_x\times\R^n_\xi).$

  (4)[Childs] difference conditions to be found in
    \cite{RS06a}.

  (5)[Coifman-Meyer] $\partial_x^\alpha\partial_\xi^\beta a(x,\xi)\in
      L^\infty(\R^n_x\times\R^n_\xi)$,
      $|\alpha|\leq [n/2]+1, \beta\in\{0,1\}^n.$

  (6)[Coifman-Meyer] $\exists 2\leq p<\infty$:
    $\partial_x^\alpha\partial_\xi^\beta a(x,\xi)\in
      L^p(\R^n_x\times\R^n_\xi)$, 
      $|\alpha|\leq [n(1/2-1/p)]+1, |\beta|\leq 2n.$

  Then $T$ is $L^2(\Rn)$--bounded.
 \end{thm}
In brackets we list the names of the authors of the corresponding
results for pseudo-differential operators.
This theorem follows from a more general statement for
the $L^2$--boundedness of Fourier integral operators with
symbols in Besov spaces that appeared in \cite{RS06a}.

There are also global $L^2$--boundedness theorems for adjoint
operators, as well as for Fourier integral operators
with general amplitudes.
We refer to \cite{RS06a} for details of all such statements,
as well as for results in weighted $L^2(\Rn)$ spaces.

The global calculus of operators \eqref{EQ:globalT} as well
as global boundedness theorems in weighted Sobolev spaces
can be found in \cite{RS06b}.

\section{Fourier integral operators in Colombeau's spaces}

In this section we will discuss properties of Fourier integral 
operators in the Colombeau's spaces of new generalised
functions. For details on these spaces we refer to e.g. 
\cite{Co84} or \cite{NPS98}. 

The theory of pseudo-differential operators in Colombeau's
spaces has been developed in \cite{GGO03} and 
\cite{GH05}. Elements of
the corresponding theory of Fourier integral operators has
been laid down in \cite{GHO07}. 

At the same time, hyperbolic partial differential equations
in Colombeau's spaces have been studied in
\cite{LO91}, \cite{Ho04} by the energy methods, yielding
relevant extensions of $L^2$ type results to the setting of
new generalised functions.

In this section we will present corresponding $L^p$--results 
in the setting of Colombeau's generalised functions,
for Fourier integral operators, with subsequent corresponding
implications for solutions to hyperbolic equations.

Let $X$ be a bounded open set. 
Consider families $(u_\epsilon)_\epsilon$ 
of functions $u_\epsilon\in L^p_\infty(X)$, $0<\epsilon\leq 1$,
where $L^p_\infty(X)$ stands for the space of function in
$L^p(X)$ for which all derivatives also belong to $L^p(X)$.
One can single out several important families of such functions
dependent on their behaviour with respect to $\epsilon$.
Thus, the class of moderate families ${\mathcal E}_{L^p}(X)$ 
is defined as the collection of families satisfying 
$$\forall \alpha\geq 0\; \exists N\geq 0:\;
||\partial^\alpha u_\epsilon||_{L^p}=O(\epsilon^{-N})
\textrm{ as } \epsilon\to 0.$$ 
The class of null families ${\mathcal N}_{L^p}(X)$
is defined as a subclass of ${\mathcal E}_{L^p}(X)$
satisfying the condition that
$$\forall N\geq 0:\;
||u_\epsilon||_{L^p}=O(\epsilon^{N})
\textrm{ as } \epsilon\to 0.$$ 
The Colombeau's algebra ${\mathcal G}_{L^p}(X)$ is then defined
as
$${\mathcal G}_{L^p}(X):={\mathcal E}_{L^p}(X)/
{\mathcal N}_{L^p}(X).$$
Distributions $L^{p}_{-\infty}=\bigcup_{s\in\R} L^{p}_{s}$
are embedded in ${\mathcal G}_{L^p}(X)$ by the mapping
$\iota(u)=[(u*(\rho_\epsilon)]_\epsilon$, where
$\rho_\epsilon(x)=\epsilon^{-n}\rho(x/\epsilon)$ is the
standard Friedrichs mollifier. 

Subsequently, one can define the corresponding
classes of generalized pseudo-differential operators.
Thus, the generalized symbol is defined as a family 
$(a_\epsilon)_\epsilon$ of usual
symbols $a_\epsilon\in S^m$ such that
$$\forall k,l\;\; \exists N:
\sup_{|\alpha|\leq k, |\beta|\leq l}
\sup_{x,\xi\in\Rn} (1+|\xi|)^{-m+|\alpha|}
\left|\partial_\xi^\alpha\partial_x^\beta a_\epsilon(x,\xi)
\right|=O(\epsilon^{-N}) \textrm{ as } \epsilon\to 0.$$

Then $A:{\mathcal G}_{L^p}(X)\to {\mathcal G}_{L^p}(X)$ is a
generalised pseudo-differential operator with generalised symbol
$(a_\epsilon)_\epsilon$ if, on the representative level,
it acts as $$(u_\epsilon)_\epsilon\mapsto
(a_\epsilon(t,x,D_x)u_\epsilon)_\epsilon.$$

To formulate the results, also the notion of the slow scale
is required. A generalised symbol $(a_\epsilon)_\epsilon$ is 
said to be of
log--type up to order $(k,l)$ (or to be slow
scale) if
$$\sup_{|\alpha|\leq k, |\beta|\leq l}
\sup_{x,\xi\in\Rn} (1+|\xi|)^{-m+|\alpha|}
\left|\partial_\xi^\alpha\partial_x^\beta a_\epsilon(x,\xi)
\right|=O(\log(1/\epsilon)) \textrm{ as } \epsilon\to 0.$$

In \cite{LO91}, Lafon and Oberguggenberger 
considered partial differential operators 
of the form
$$A=\sum_{j=1}^n a_j(t,x)\partial_{x_j}+b(t,x).$$ 
They investigated the Cauchy problem for operator $D_t+A$
and showed the existence of solutions 
in ${\mathcal G}_{L^\infty}(X)$ provided that  
$b$ and $\partial_{x_k}a_j$ are of log--type.
Moreover, if $a_j$ and $b$ are constant for large $x$, then
the solution is also unique, and in 
\cite{Ob89} an example of the non--uniqueness was given
if this condition breaks. 
In \cite{Ho04}, H\"ormann considered the Cauchy problem
for more general pseudo-differential operators $A$
of log--type, for which he showed 
existence and uniqueness in ${\mathcal G}_{L^2}(X)$
(also giving some estimates on $k$ and $l$ in the log-type 
assumption).
In particular, the non--uniqueness effect disappears in 
${\mathcal G}_{L^2}(X)$ compared to ${\mathcal G}_{L^\infty}(X)$.

These results may provide some hints on the behaviour of 
relevant Fourier integral operators in Colombeau' spaces.
Let us consider generalised Fourier integral operators
of the form 
$$Tu(x)=\int_{X}\int_{\R^n} e^{i\phi(x,y,\xi)}
a(x,y,\xi) u(y) dy \dslash\xi.$$
Let us define the regular Colombeau's algebra.
The class of regular families ${\mathcal R}_{L^p}(X)$ is defined
as a class of functions satisfying the condition
$$\exists N\geq 0 \;\forall \alpha\geq 0:\;
||\partial^\alpha u_\epsilon||_{L^p}=O(\epsilon^{-N})
\textrm{ as } \epsilon\to 0.$$ 
Then the regular Colombeau algebra ${\mathcal G}^\infty_{L^p}(X)$
is defined by
$${\mathcal G}^\infty_{L^p}(X):=
{\mathcal R}_{L^p}(X)/{\mathcal N}_{L^p}(X).$$
In \cite{GHO07}, Garetto, H\"ormann and Oberguggenberger 
showed that if $\phi$ is 
a non--degenerate phase function then $T$ maps
${\mathcal G}_{L^\infty}(X)$ to itself continuously. 
If the phase function  $\phi$ is a generalised family
satisfying the slow scale assumption (e.g. log--type),
and $a=(a_\epsilon)_\epsilon$ is a regular family of
amplitudes, then $T$ maps locally continuously 
the space of regular Colombeau's 
functions ${\mathcal G}^\infty_{L^\infty}(X)$ to itself.
The following theorem extends this to Colombeau's algebras
over any $L^p(X)$, $1\leq p\leq\infty$.

\begin{thm} Let $1\leq p\leq\infty$.
If $\phi$ is 
a non--degenerate (generalised) phase function then $T$ maps
${\mathcal G}_{L^p}(X)$ to itself continuously. 
If the phase function  $\phi$ is a generalised family
satisfying the slow scale assumption (e.g. log--type),
and $a=(a_\epsilon)_\epsilon$ is a regular family of
amplitudes, then $A$ maps locally continuously 
the space of regular Colombeau's 
generalised functions ${\mathcal G}^\infty_{L^p}(X)$ to itself.
\end{thm}
Consequently, propagators for strictly
hyperbolic Cauchy problems are continuous in
${\mathcal G}_{L^p}(X)$, and in ${\mathcal G}^\infty_{L^p}(X)$
under the log--type assumption on the phase, for all 
$1\leq p\leq\infty$.

The proof is based on the extension to the Colombeau's setting 
of the eikonal and transport equations, modulo controllable
errors with respect to $\epsilon$, and on estimates for 
generalised Fourier integral operators in 
${\mathcal G}_{L^p}(X)$. Details of proofs and exact losses
of regularity in Colombeau's spaces will appear
in \cite{Rup}.


\end{document}